\documentstyle[12pt]{article}
\newcommand{\sect}[1]{\section{#1}\setcounter{equation}{0}}

\font\mbn=msbm10 scaled \magstep1
\font\mbs=msbm7 scaled \magstep1
\font\mbss=msbm5 scaled \magstep1
\newfam\mbff
\textfont\mbff=\mbn
\scriptfont\mbff=\mbs
\scriptscriptfont\mbff=\mbss\def\mbf{\fam\mbff}
\def\Re{{\mbf R}}

\def\Z{{\mbf Z}}
\def\Co{{\mbf C}}

\def\P{{\mbf P}}

\def\N{{\mbf N}}

\newtheorem{Th}{Theorem}[section]
\newtheorem{Lm}[Th]{Lemma}
\newtheorem{C}[Th]{Corollary}
\newtheorem{D}[Th]{Definition}
\newtheorem{Proposition}[Th]{Proposition}
\newtheorem{R}[Th]{Remark}

\author{Alexander Brudnyi\thanks{Research supported in part by NSERC.
}\\
Department of Mathematics and Statistics\\
University of Calgary, Calgary\\
Canada}
\title{On Center Sets of ODEs Determined by Moments of their 
Coefficients}
\date{} 
\begin{document} 
\maketitle
\begin{abstract}
{The classical H. Poincar\'{e} Center-Focus problem asks about 
the characterization of planar polynomial vector fields such that all their
integral trajectories are closed curves whose interiors contain a fixed point,
a {\em center}. This problem can be reduced to a center problem for
some ordinary differential equation whose coefficients are trigonometric
polynomials depending polynomially on the coefficients of the field. In
this paper we show that the set of centers in the Center-Focus problem
can be determined as the set of zeros of some continuous functions from the
moments of coefficients of this equation.}
\end{abstract}

\sect{\hspace*{-1em}. Introduction.}
Consider the ordinary differential equation 
\begin{equation}\label{e1}
\frac{dv}{dx}=\sum_{i=1}^{\infty}a_{i}(x)\ \!v^{i+1}\ ,\ \ \ \ \
x\in I_{T}:=[0,T]\ ,
\end{equation}
whose coefficients $a_{i}$ belong to the Banach space $L^{\infty}(I_{T})$ of 
bounded measurable complex-valued functions on $I_{T}$ with the supremum norm.
Let $X_{i}:=L^{\infty}(I_{T})$ be the space of coefficients $a_{i}$ from
(\ref{e1}), and $X$ be the complex Fr\'{e}chet space of sequences
$a=(a_{1},a_{2},\dots)\in\prod_{i\geq 1}X_{i}$ satisfying 
\begin{equation}\label{e11}
\sup_{x\in I_{T}, i\in\N}\sqrt[i]{|a_{i}(x)|}<\infty.
\end{equation}
Condition (\ref{e11}) implies that for any $a\in X$ the corresponding
equation (\ref{e1}) has Lipschitz solutions on $I_{T}$ for all sufficiently 
small initial values. Such solutions can be obtained by Picard iteration.
We say that equation (\ref{e1}) determines a {\em center} if every solution
$v$ of (\ref{e1}) with a sufficiently small initial value satisfies 
$v(T)=v(0)$.
By ${\cal C}\subset X$ we denote the set of centers. 
The center problem for equation (\ref{e1}) is: 
{\em given $a\in X$ to determine
whether $a\in {\cal C}$}.
This problem arises naturally in the framework of the geometric theory
of ordinary differential equations created by Poincar\'{e}.
In particular, there is a connection between this problem and the classical 
Poincar\'{e} Center-Focus problem for polynomial vector fields. 
For a non-degenerate equilibrium point of such a field it can
be reduced to the following equivalent form.

Consider the system of ODEs in the plane
\begin{equation}\label{e2}
\frac{dx}{dt}=-y+F(x,y)\ ,\ \ \ \frac{dy}{dt}=x+G(x,y)
\end{equation}
where $F,G$ are real analytic functions in a neigbourhood of $0\in\Re^{2}$
whose Taylor expansions at $0$ do not contain constant and linear terms.
Then for $F,G$ polynomials of a given degree, the Poincar\'{e}
Center-Focus Problem asks about conditions on the coefficients of $F$ and
$G$ under which all trajectories of (\ref{e2}) situated in a small 
neighbourhood of the origin are closed. (A similar problem can be posed
for the general case.) Observe that passing to polar coordinates
$x=r\cos\phi,\ y=r\sin\phi$ in (\ref{e2}) and expanding the right-hand
side of the resulting equation as a series in $r$ (for
$F$ and $G$ whose coefficients are sufficiently small) we obtain an 
equation of the form (\ref{e1}) whose coefficients are trigonometric 
polynomials
depending polynomially on the coefficients of (\ref{e2}).
This reduces the Center-Focus Problem for (\ref{e2}) to the center problem for
equations (\ref{e1}) with coefficients depending polynomially on a
parameter. In this paper we will consider a general equation of such form
\begin{equation}\label{e3}
\frac{dv}{dx}=\sum_{i=1}^{\infty}a_{i}(w,x)\ \! v^{i+1}\ ,\ \ \ \ 
x\in I_{T}\ ,
\end{equation}
where all $a_{i}$ are holomorphic functions in $w$ varying in a Stein domain 
$U\subset\Co^{k}$ and $a_{i}(w,\cdot)\in L^{\infty}(I_{T})$. 
We set
$A(w):=(a_{1}(w,\cdot),a_{2}(w,\cdot),\dots)$ and assume that 
$A(w)\in X$ for all $w\in U$. This gives a holomorphic map 
$A: U\to X$. Our goal is to study the center problem for
(\ref{e3}), i.e., the problem of the characterization of the set 
${\cal C}_{A}:=A^{-1}({\cal C})\subset U$.
\sect{\hspace*{-1em}. Formulation of Main Results.}
{\bf 2.1. A Characterization of the Center Set.} In [Br1] a 
characterization of the center set ${\cal C}$ is
given in terms of the iterated path integrals. For $a\in X$ let us 
consider the {\em basic iterated integrals}
\begin{equation}\label{iter}
I_{i_{1},\dots,i_{k}}(a):=
\int\cdots\int_{0\leq s_{1}\leq\cdots\leq s_{k}\leq T}
a_{i_{k}}(s_{k})\cdots a_{i_{1}}(s_{1})\ ds_{k}\cdots ds_{1}
\end{equation}
(for $k=0$ we assume that this equals 1).
They can be thought of as $k$-linear holomorphic functions on $X$. By
the Ree formula [R] the linear space generated by
all such functions is an algebra. A linear combination of iterated
integrals of order $\leq k$  is called an {\em iterated polynomial of degree 
$k$}. 

Let $v(x;r;a)$, $x\in I_{T}$, be the solution of equation 
(\ref{e1}) corresponding to $a\in X$ with (sufficiently small)
initial value $v(0;r;a)=r$. Then $P(a)(r):=v(T;r;a)$ is the
{\em first return map} of this equation.\\
\\
{\bf Theorem ([Br1])} {\em
For sufficiently small initial values $r$ the first return map $P(a)$ is an
absolutely convergent power series 
$P(a)(r)=r+\sum_{n=1}^{\infty}c_{n}(a)r^{n+1}$,
where
$$
c_{n}(a)=\sum_{i_{1}+\dots + i_{k}=n}c_{i_{1},\dots,i_{k}}
I_{i_{1},\dots,i_{k}}(a)\ , \ \ \ {\rm and}
$$
$$
c_{i_{1},\dots,i_{k}}=
(n-i_{1}+1)\cdot (n-i_{1}-i_{2}+1)\cdot (n-i_{1}-i_{2}-i_{3}+1)\cdots 1\ .
$$

The center set ${\cal C}\subset X$ of equation (\ref{e2}) is determined
by the system of polynomial equations $c_{n}(a)=0$, $n=1,2,\dots$ .}

This theorem implies that {\em the center set ${\cal C}_{A}$ for
(\ref{e3}) is a closed complex subvariety of $U$ determined as
the set of zeros of holomorphic functions
$A^{*}c_{n}$, $1\leq n\leq N$;
here $A^{*}c_{n}(w):=c_{n}(A(w))$, $w\in U$.}\\
{\bf 2.2. Algebraic Structure of the Center Set.} In [Br3] we constructed
an algebraic model for the center problem for equation (\ref{e1}). Let us
briefly describe this construction.

First, from the Theorem of the previous section it follows that any point
$a\in X$ such that all iterated integrals vanish at $a$, belongs to 
${\cal C}$.
The set ${\cal U}\subset {\cal C}$ of all such points is called the
{\em set of universal centers} of (\ref{e1}). In [Br2] we characterized
some elements from ${\cal U}$. It was related to a certain composition
condition whose role and importance for Abel differential equations
with polynomial coefficients was studied in [AL], [BFY1], [BFY2], [Y].

Next, given $a=(a_{1},a_{2},\dots)$ and $b=(b_{1},b_{2},\dots)$ 
from $X$ we define 
$$
a*b=(a_{1}*b_{1},a_{2}*b_{2},\dots)\in X\ \ \
{\rm and}\ \ \ a^{-1}=(a_{1}^{-1},a_{2}^{-1},\dots)\in X
$$
where for $i\in\N$,
$$
(a_{i}*b_{i})(t)=\left\{
\begin{array}{ccc}
2a_{i}(2t)&{\rm if}&0<t\leq T/2\\
2b_{i}(2t-T)&{\rm if}&T/2<t\leq T
\end{array}
\right.
$$
and
$$
a_{i}^{-1}(t)=-a_{i}(T-t)\ ,\ \ \ 0<t\leq T\ .
$$

We say that $a,b\in X$ are
{\em equivalent} (written, $a\sim b$) if $a*b^{-1}\in {\cal U}$. 
It was shown in [Br3]
that $\sim$ is an equivalence 
relation, that is  $X$ partitions into mutually disjoint equivalence classes. 
Let $G(X)$ be the set of these classes.
Then $*$ induces a 
multiplication $\cdot :G(X)\times G(X)\rightarrow G(X)$ converting
$G(X)$ to a group. Also, 
the iterated integrals are constant on any equivalence class. 
So they can be considered as functions on 
$G(X)$. It was proved that these functions separate points on $G(X)$. 
Let us equip $G(X)$ with the weakest topology in which all 
iterated integrals are continuous. 
Then $G(X)$ is a separable topological group. Moreover,
it is contractible, arcwise connected, locally arcwise and simply connected
and residually torsion free nilpotent (i.e., 
finite-dimensional unipotent representations separate elements of $G(X)$).

Further, consider the set $G_{c}[[r]]$ of complex power series 
$f(r)=r+\sum_{i=1}^{\infty}d_{i}r^{i+1}$ each convergent in some  
neighbourhood of $0\in\Co$. Let $d_{i}:G_{c}[[r]]\rightarrow\Co$ be such 
that $d_{i}(f)$ equals the $(i+1)$-st coefficient of the Taylor expansion of 
$f$ at $0$. We equip $G_{c}[[r]]$ with the weakest topology in which all 
$d_{i}$ are continuous functions and consider the multiplication $\circ$ 
on $G_{c}[[r]]$ defined by the composition of series. Then 
$G_{c}[[r]]$ is a separable topological group. Moreover, it is contractible,
arcwise connected, locally arcwise and simply connected and
residually torsion free nilpotent.

Now, for any $a\in X$ let $v(x;r;a)$, $x\in I_{T}$, be the
Lipschitz solution of equation (\ref{e2}) corresponding to $a$ with 
initial value $v(0;r;a)=r$. Clearly for every $x$ we have 
$v(x;r;a)\in G_{c}[[r]]$. Let $P(a)(r):=v(T;r;a)$ be the first return
map. Then 
$$
P(a*b)=P(b)\circ P(a)\ .
$$
This together with the fact that $P(a)(r)\equiv r$ for 
any $a\sim 0\ (\in X)$ imply
that there is a map $\widehat P:G(X)\rightarrow G_{c}[[r]]$
such that $\widehat P([a]):=P(a)$ where $[a]$ denotes the equivalence class
containing $a\in X$. It was proved that
$\widehat P$ is a surjective homomorphism of topological groups,
the kernel $\widehat {\cal C}\subset G(X)$ of $\widehat P$
coincides with the image of the center set ${\cal C}\subset X$ in $G(X)$,
and ${\cal C}$ is contractible,
arcwise connected,  locally arcwise and simply connected. Also, there is a 
continuous map $T:G_{c}[[r]]\rightarrow G(X)$ such that 
$\widehat P\circ T=id$. In particular, the map 
$\widetilde T:G_{c}[[r]]\times\widehat {\cal C}\rightarrow G(X)$,
$\widetilde T(f,g):=T(f)\cdot g$, is a homeomorphism.\\
{\bf 2.3. Moments.}
It is natural to try to characterize the center set ${\cal C}_{A}$ of
(\ref{e3}) as the
set of zeros of functions from iterated integrals of a simple form. One of
the simplest type of 
iterated integrals are called {\em moments}. Moments play
an important role in the study of the center problem for Abel differential
equations (see, e.g., [AL], [BFY1], [BFY2], [Y]). 
For $a=(a_{1},a_{2},\dots)\in X$ they are defined by the formulas
\begin{equation}\label{moment}
m_{i_{1},\dots,i_{k+1}}^{n_{1},\dots,n_{k}}(a):=
\int_{0}^{T}(\widetilde a_{i_{1}}(s))^{n_{1}}\cdots
(\widetilde a_{i_{k}}(s))^{n_{k}}a_{i_{k+1}}(s)\ \! ds
\end{equation}
for all possible natural numbers $n_{1},\dots, n_{k}$ and 
$i_{1},\dots, i_{k+1}$ where
\begin{equation}\label{moment1}
\widetilde a_{i}(x):=\int_{0}^{x}a_{i}(s)\ \!ds\ .
\end{equation}
By ${\cal M}$ we denote the set of moments on $X$.
One can consider elements of ${\cal M}$ as continuous functions on
$G(X)$ (see section 2.2). By $X_{*}\subset X$ we denote the set of elements $a\in X$ such that $I_{k}(a)=0$ for all $k\in\N$. Let $G(X_{*})$ be the image of $X_{*}$ in $G(X)$. Then $G(X_{*})$ is a normal subroup of $G(X)$.
It is easy to see that every $m\in {\cal M}$
satisfies $m(g_{1}g_{2})=m(g_{1})+m(g_{2})$
for all $g_{1},g_{2}\in G(X_{*})$. Thus ${\cal M}$ consists of continuous 
characters of $G(X_{*})$. 
Set $G^{1}(X_{*}):=[G(X_{*}),G(X_{*})]$ and let $\overline{G^{1}(X_{*})}$ be the closure of
$G^{1}(X_{*})$ in $G(X)$. In a forthcoming paper we will show that 
$g\in\overline{G^{1}(X_{*})}$ if and only if $m(g)=0$ for any $m\in {\cal M}$.
In particular, elements of ${\cal M}$ separates points on 
$G(X_{*})/\overline{G^{1}(X_{*})}$. However the 
difference between $G^{1}(X_{*})$ and $\overline{G^{1}(X_{*})}$ is not yet well 
understood. We will just mention the following important
problem of independent interest: {\em to describe points $g\in G(X)$
such that if $m(g)=0$ for all $m\in {\cal M}$ then $g\in G^{1}(X_{*})$.}
Some results in this direction can be obtained from the description of
polynomially convex hulls of smooth curves in $\Co^{n}$ (see [W], [A],
[HL] and section 3.1 below).\\
{\bf 2.4. Main Results.} We retain the notation of the Introduction.
Let us consider the equation (\ref{e3}) with all coefficients being either
polynomials in $x$ or in $e^{\pm 2\pi ix/T}$. Then we prove
\begin{Th}\label{te1}
For every compact set $K\subset U$ there are a number $N=N(K)\in\N$,
moments $m_{1},\dots, m_{N}\in {\cal M}$ and continuous functions
$f_{1},\dots,f_{N}$ on $\Co^{N}$ such that the center set
${\cal C}_{A}\cap K$ of (\ref{e3}) in $K$ is
determined as the set of zeros of continuous functions 
$f_{n}(A^{*}m_{1},\dots, A^{*}m_{N})$ on $U$ $(1\leq n\leq N)$.
Moreover, the set of universal centers 
$A^{-1}({\cal U})\cap K$ of (\ref{e3}) in $K$
is the set of zeros of the holomorphic functions
$A^{*}m_{n}$, $1\leq n\leq N$, restricted to $K$.
\end{Th}

To formulate some specific versions of this result we need the following
\begin{D}\label{d0}
A holomorphic polynomial $p$ on $\Co^{k}$ is called 
{\em $\alpha$-homogeneous} for some 
$\alpha=(\alpha_{1},\dots,\alpha_{k})\in\N^{k}$ if 
$\widetilde p(w_{1},\dots,w_{k}):=p(w_{1}^{\alpha_{1}},\dots 
w_{k}^{\alpha_{k}})$, $(w_{1},\dots,w_{k})\in\Co^{k}$, is a homogeneous
polynomial.
\end{D}

Assume that all coefficients $a_{i}$ in (\ref{e3}) are $\alpha$-homogeneous
polynomials in $w\in\Co^{k}$ for some $\alpha\in\N^{k}$ whose coefficients
are either polynomials in $x$ or in $e^{\pm 2\pi x/T}$. Assume also
that the set of universal centers 
${\cal U}_{A}:=A^{-1}({\cal U})\subset\Co^{k}$ of (\ref{e3})
coincides with $0\in\Co^{k}$ (this holds, e.g., if there are 
$a_{i_{1}},\dots, a_{i_{l}}$ such that the set of zeros of 
$\widehat a_{i_{s}}(w):=\int_{0}^{T}a_{i_{s}}(w,x)\ \!dx$, $1\leq s\leq l$,
$w\in\Co^{k}$, coincides with $0$.) Then we prove
\begin{Th}\label{cor1}
There are moments $m_{1},\dots,m_{N}\in {\cal M}$ and holomorphic
polynomials $p_{1},\dots, p_{N}$ on $\Co^{N}$ such that the center set
${\cal C}_{A}$ of (\ref{e3}) is determined as the set of zeros of holomorphic
polynomials $p_{n}(A^{*}m_{1},\dots,A^{*}m_{N})$, $1\leq n\leq N$.
\end{Th}

Assume as before that all coefficients $a_{i}$ in (\ref{e3}) are 
$\alpha$-homogeneous polynomials in $w\in\Co^{k}$ for some 
$\alpha\in\N^{k}$ whose coefficients are either polynomials in $x$ or in 
$e^{\pm 2\pi i x/T}$. Then the set of universal centers 
${\cal U}_{A}$ of (\ref{e3}) is a complex algebraic subvariety of $\Co^{k}$. 
By $\overline {\cal U}_{A}$ we denote its closure in 
$\Co\P^{k}:=\Co^{k}\cup H_{k}$, the complex 
$k$-dimensional projective space (here $H_{k}$ is the complex hyperplane at 
infinity). Then $\overline {\cal U}_{A}\subset\Co\P^{k}$ is a projective 
subvariety. Let $V\subset\Co^{k}$ be
a complex algebraic subvariety whose closure in $\Co\P^{k}$ does not intersect
$\overline {\cal U}_{A}$. 
\begin{Th}\label{cor2}
There are moments $m_{1},\dots,m_{N}\in {\cal M}$ and holomorphic
polynomials $p_{1},\dots, p_{N}$ on $\Co^{N}$ such that the center set
${\cal C}_{A}\cap V$ of (\ref{e3}) in $V$ is the set of
zeros of holomorphic polynomials $p_{n}(A^{*}m_{1},\dots,A^{*}m_{N})$, 
$1\leq n\leq N$, restricted to $V$.
\end{Th}
\begin{R}\label{re0}
{\rm (1)  All coefficients in (\ref{e3}) 
obtained from complex polynomial vector fields (\ref{e2}) of degree $d$ are 
$\alpha$-homogeneous polynomials on $\Co^{k}$; here $k:=d^{2}+3d-4$ is the 
number of coefficients of polynomials $F$ and $G$ in (\ref{e2}).\\
(2) For Abel differential equations with polynomial 
coefficients the result of Theorem \ref{cor1} was conjectured by Yomdin 
[Y1] without any additional assumption on ${\cal U}_{A}$.\\
(3) Under the hypotheses of Theorem \ref{te1} one can give a description of 
elements of ${\cal U}_{A}$ in terms of the composition condition, see, e.g., 
[Br2, Corollaries 1.19, 1.20].\\
(4) Since $\widehat {\cal C}\cap \overline{G^{1}(X)}\neq\{1\}$,
Theorem \ref{te1} does not hold for a general equation of the form 
(\ref{e3}).}
\end{R}
\sect{\hspace*{-1em}. Separation by Moments.}
{\bf 3.1. Auxiliary Results.} 
We require
\begin{D}\label{de1}
The polynomially convex hull $\widehat K$ of a compact set
$K\subset\Co^{n}$ is the set of points $z\in\Co^{n}$
such that if $p$ is any holomorphic polynomial in $n$ variables
$$
|p(z)|\leq\max_{x\in K}|p(x)|\ .
$$
\end{D}
It is well known (see e.g. [AW]) that $\widehat K$ is compact, and if
$K$ is connected then $\widehat K$ is connected.

Let $\widetilde A=(\widetilde a_{1},\dots,\widetilde a_{n}): 
I_{T}\to\Co^{n}$ be a Lipschitz curve where $\widetilde a_{i}$ are the first
integrals of some functions $a_{i}\in L^{\infty}(I_{T})$, 
see (\ref{moment1}). Set $\Gamma:=\widetilde A(I_{T})$. Then arguments
similar (but essentially more easier) to those used in [Br2, Theorem 1.10]
give
\begin{Th}\label{te11}
Assume that all moments from functions $a_{1},\dots, a_{n}$ are zeros.
Then the path $\widetilde A: I_{T}\to\Co^{n}$ is closed and represents 0
in the homology groups $H_{1}(U,\Co)$ for any connected neighbourhood
$U$ of $\widehat\Gamma$.\ \ \ \ \ $\Box$
\end{Th}
Theorem \ref{te11} can be sharpened under some additional assumptions on 
$\widehat\Gamma$. For instance, similarly to the proof of
[Br2, Corollary 1.12] one obtains
\begin{C}\label{co1}
Assume that $\widetilde A: I_{T}\to\Co^{n}$ satisfies hypotheses of
Theorem \ref{te11} and
$\widehat\Gamma$ is triangulable. Then 
$\widetilde A: I_{T}\to\widehat\Gamma$ represents 0 in 
$H_{1}(\widehat\Gamma,\Co)$.\ \ \ \ \ $\Box$
\end{C}
From the result of Alexander [A] it follows that
the covering dimension of $\widehat\Gamma$ is 2. However, it does not imply
its triangulability. Some cases when $\widehat\Gamma$ is triangulable are
considered in [W], [HL] (see also [Br2]).\\
{\bf 3.2.} The proof of Theorem \ref{te1} is based on a 
separation result proved in this section.

Let ${\cal E}_{n}$ be the set of closed complex 
curves in $\Co^{n}$ containing 0  whose normalizations are either
$\Co$ or $\Co^{*}\ (:=\Co\setminus\{0\})$. Consider the equation
\begin{equation}\label{e5}
\frac{dv}{dx}=\sum_{i=1}^{n}a_{i}(x)\ \!v^{i+1}\ ,\ \ \ \ x\in I_{T}\ ,\ 
a_{i}\in L^{\infty}(I_{T})\ ,\ 1\leq i\leq n\ .
\end{equation}
For $a=(a_{1},\dots, a_{n},0,\dots)\in X$ we define the map 
$\widetilde A_{a}:I_{T}\to\Co^{n}$ by the formula
\begin{equation}\label{e6}
\widetilde A_{a}(x):=(\widetilde a_{1}(x),\dots,\widetilde a_{n}(x))
\end{equation}
with $\widetilde a_{i}$ given by (\ref{moment1}). We will assume that
\begin{itemize}
\item[$(C_{n1})$]
There is a curve $C\in {\cal E}_{n}$ such that 
$\widetilde A_{a}(I_{T})\subset C$.
\item[$(C_{n2})$]
There is a continuous map $A_{a}':I_{T}\to C_{n}$ such that 
$\widetilde A_{a}=n\circ A_{a}'$ where $n:C_{n}\to C$ is the normalization 
of $C$.
\end{itemize}
\begin{R}\label{r1}
{\rm Condition $(C_{n2})$ is valid, for instance, if $\widetilde A_{a}$ can be
extended to a holomorphic map $U_{a}\to C$ where
$U_{a}$ is a neighbourhood of $I_{T}\subset\Co$.}
\end{R}

Let $D_{n}$ denote the class of equations (\ref{e5}) satisfying 
conditions $(C_{n1})$ and $(C_{n2})$. 
By $X(D_{n})\subset X$ we denote the set of
sequences $a=(a_{1},\dots, a_{n},0,\dots)\in X$ whose coordinates
are coefficients of equations from $D_{n}$. 
As before $[a]\in G(X)$ stands for the equivalence
class containing $a\in X$.
\begin{Th}\label{te2}
Let $a,b\in X(D_{n})$ be such that $m([a])=m([b])$ for every moment 
$m\in {\cal M}$. Then $[a]=[b]$.
\end{Th}
\begin{R}\label{r2}
{\rm
Theorem \ref{te2} implies that if $a\in X(D_{n})$ is such that
$m(a)=0$ for every $m\in {\cal M}$, then $a\in {\cal U}$, the set of universal
centers. In this case the structure of $a$ can be described explicitly 
(see [Br2, Corollary 1.17]).}
\end{R}
Since every moment vanishes on the
commutant $G^{1}(X)$, Theorem \ref{te2} implies
\begin{C}\label{c1}
Let $a,b\in X(D_{n})$ satisfy $[a*b^{-1}]\in G^{1}(X)$. Then
$[a]=[b]$.\ \ \ \ \ $\Box$
\end{C}
Since the center set ${\cal C}$ contains ${\cal U}$ and is closed
under the multiplication $*$, we get
\begin{C}\label{c2}
Let $a,b\in X(D_{n})$ be such that $a\in {\cal C}$ and $m([a])=m([b])$
for every moment $m\in {\cal M}$. Then $b\in {\cal C}$.\ \ \ \ \
$\Box$
\end{C}
{\bf Proof of Theorem \ref{te2}.}
For $a$ and $b$ as in Theorem \ref{te2} we consider
the map $\widetilde A_{a*b^{-1}}:I_{T}\to\Co^{n}$ defined by (\ref{e6}). Then
the hypotheses of the theorem imply that
\begin{itemize}
\item[(1)]
$$
\widetilde A_{a*b^{-1}}(0)=\widetilde A_{a*b^{-1}}(T)=0\ ;
$$
\item[(2)]
There are closed complex curves $C_{1}$ and $C_{2}$ from ${\cal E}_{n}$ 
such that 
$$
\widetilde A_{a*b^{-1}}(I_{T})\subset C_{1}\cup C_{2}\ \ \ {\rm and}\ \ \
\widetilde A_{a*b^{-1}}(I_{l})\subset C_{l}\ \! ,\ \ \ l=1,2\ ,
$$
where $I_{1}:=[0,T/2]$ and $I_{2}:=[T/2,T]$.
\end{itemize}
Extending $\widetilde A_{a*b^{-1}}$ to a $T$-periodic map $\Re\to\Co^{n}$ 
and identifying the quotient $\Re/(T\cdot\Z)$ with a circle 
$S$ we think of $\widetilde A_{a*b^{-1}}$ as a closed path $S\to\Co^{n}$.
Set $\Gamma:=\widetilde A_{a*b^{-1}}(I_{T})$ and 
$\Gamma_{l}:=\widetilde A_{a*b^{-1}}(I_{l})$, $l=1,2$.
Let $\widehat\Gamma\subset\Co^{n}$ be the polynomially convex hull of 
$\Gamma$. Then from Theorem \ref{te11} by the hypotheses of Theorem \ref{te2}
we get that {\em for any connected neighbourhood $U$ of
$\widehat\Gamma$ the path
$\widetilde A_{a*b^{-1}}:S\to U$ represents 0 in $H_{1}(U,\Co)$}. Next,
since $C_{1}\cup C_{2}$ is a closed complex subvariety of $\Co^{n}$, an 
argument used in [Br2, Lemma 5.1] shows that
$\widehat\Gamma\subset C_{1}\cup C_{2}$.  Since $C_{1}\cup C_{2}$ is
triangulable, the previous implies easily that
$\widetilde A_{a*b^{-1}}:S\to C_{1}\cup C_{2}$ represents 0 in 
$H_{1}(C_{1}\cup C_{2},\Co)$. We will consider two cases.
\begin{itemize}
\item[{\bf (A)}]\ \ \
$C_{1}\neq C_{2}$.
\end{itemize}

Let us prove
\begin{Lm}\label{le1}
$$
\widetilde A_{a*b^{-1}}(T/2)=0\ .
$$
\end{Lm}
{\bf Proof.} Suppose, to the contrary, that 
$\widetilde A_{a*b^{-1}}(T/2):=z\neq 0$.
Let $f$ be a function on $C_{1}$ constant in some
neighbourhoods of singular points of $C_{1}\cup C_{2}$ in $C_{1}$, 
smooth outside these singular points, equals 0 in a neigbourhood of 
$0\in\Co^{n}$ and 1 in a neighbourhood of $z$. 
Then $df$ is well defined smooth
1-form on $C_{1}$. We can extend it by 0 to $C_{2}$ to get a smooth
$d$-closed 1-form $\omega$ on $C_{1}\cup C_{2}$. Since $C_{1}\cup C_{2}$ is
a complex space, it is Lipschitz triangulable. This and the fact that
$\widetilde A_{a*b^{-1}}:S\to C_{1}\cup C_{2}$ represents 0 in 
$H_{1}(C_{1}\cup C_{2},\Co)$ imply that
$\int_{0}^{T}\widetilde A_{a*b^{-1}}^{*}(\omega)=0$. On the other hand,
$$
\int_{0}^{T}\widetilde A_{a*b^{-1}}^{*}(\omega)=\int_{0}^{T/2}
\widetilde A_{a*b^{-1}}^{*}(df)=
f(z)-f(0)=1\ .
$$
This contradiction proves the lemma.\ \ \ \ \ $\Box$

From this lemma we obtain that $\widetilde A_{a*b^{-1}}:I_{l}\to C_{l}$, 
$l=1,2$, are closed paths.
\begin{Lm}\label{le2}
$\widetilde A_{a*b^{-1}}:I_{l}\to C_{l}$
represents 0 in $H_{1}(C_{l},\Co)$, $l=1,2$.
\end{Lm}
{\bf Proof.} Let $U_{l}$ be a neighbourhood of $C_{l}$ such that
$C_{l}$ is a deformation retract of $U_{l}$. (Such $U_{l}$  exists, 
e.g., by the Lojasiewicz triangulation theorem [L].) Then 
$H_{1}(U_{l},\Co)\cong H_{1}(C_{l},\Co)$. Also we can define
$H_{1}(U_{l},\Co)$ by the de Rham theorem (i.e., by integration of 
smooth $d$-closed 1-forms on $U_{l}$). Let $S_{l}\subset C_{l}$ be the
(discrete) set of singular points of $C_{1}\cup C_{2}$ in $C_{l}$. Then it is
easy to see that for every smooth $d$-closed 1-form $\omega$ on $U_{l}$
there exists a smooth $d$-closed 1-form $\omega'$ on $U_{l}$, equals 0 
in a neigbourhood of $S_{l}$, such that $\omega-\omega'=df$ for some
smooth function $f$ on $U_{l}$. Thus one can define $H_{1}(U_{l},\Co)$ by
means of such forms $\omega'$. (The class of such forms is denoted
by $E^{1}(S_{l})$.) 
Next, for $\omega\in E^{1}(S_{l})$ consider its restriction 
$\omega|_{C_{l}}$ to $C_{l}$. Then $\omega|_{C_{l}}$ can be 
extended by 0 to $C_{2}$ to get a smooth $d$-closed 1-form 
$\widetilde\omega$ on $C_{1}\cup C_{2}$. Since 
$\widetilde A_{a*b^{-1}}:S\to C_{1}\cup C_{2}$ represents 0 in 
$H_{1}(C_{1}\cup C_{2},\Co)$, as in the proof of Lemma \ref{le1} we 
obtain
$$
0=\int_{0}^{T}\widetilde A_{a*b^{-1}}^{*}(\widetilde\omega)=\int_{I_{l}}
\widetilde A_{a*b^{-1}}^{*}(\omega)\ .
$$
This completes the proof of the lemma.\ \ \ \ \ $\Box$

Note that after a reparametrization the path 
$\widetilde A_{a*b^{-1}}:I_{1}\to C_{1}$
coincides with $\widetilde A_{a}$ and the path 
$\widetilde A_{a*b^{-1}}:I_{2}\to C_{2}$ with
$-\widetilde A_{b}$.
\begin{Lm}\label{le21}
$$
[a]=[b]=1\ .
$$
\end{Lm}
{\bf Proof.}
Let  $n_{l}:C_{ln}\to C_{l}$ be the normalization of $C_{l}$. According to
condition $(C_{n2})$ in the definition of the class $X(D_{n})$ and the
above remark there are paths $A_{l}: I_{l}\to C_{ln}$ such that
$n_{l}\circ A_{l}=\widetilde A_{a*b^{-1}}|_{I_{l}}$. Since
$\widetilde A_{a*b^{-1}}:I_{l}\to C_{l}$ represents 0 in $H_{1}(C_{l},\Co)$,
the path $A_{l}$ is closed. (For otherwise, as before we can construct a
smooth $d$-closed 1-form on $C_{l}$ whose integral over 
$\widetilde A_{a*b^{-1}}|_{I_{l}}$ is not zero.) This implies easily that
$A_{l}:I_{l}\to C_{ln}$ represents 0 in $H_{1}(C_{ln},\Co)$. Then
by the definition of the class ${\cal E}_{n}$ we obtain that 
$A_{l}:I_{l}\to C_{ln}$ is homotopically trivial. Hence the same
is valid for $\widetilde A_{a}$ and $\widetilde A_{b}$. 
From here we deduce easily that all iterated 
integrals vanish at $[a]$ and $[b]$, i.e., $[a]=[b]=1$ 
(see [Br2, Corollary 1.17] for details).\ \ \ \ \ $\Box$

This completes the proof of the theorem in case (A).
\begin{itemize}
\item[{\bf (B)}]\ \ \ $C_{1}=C_{2}$.
\end{itemize}

In this case we apply condition $(C_{n2})$ and the 
arguments from the proof of Lemma \ref{le21} to get $[a*b^{-1}]=1$. 
We leave the details to the readers.

The proof of the theorem is complete.\ \ \ \ \ $\Box$
\sect{\hspace*{-1em}. Proof of Theorem \ref{te1}.}
Let $A(w):=(a_{1}(w,\cdot),\dots)\in X$, $w\in U$, be sequences 
of coefficients of equation (\ref{e3}). We set
${\cal A}:=\{[A(w)]\in G(X)\ :\ w\in U\}$. Our proof
is based on 
\begin{Proposition}\label{pr1}
Suppose that the coefficients of (\ref{e3}) satisfy the hypotheses of 
Theorem \ref{te1}. Then moments separate points on ${\cal A}$.
\end{Proposition}
{\bf Proof.} Assume to the contrary that there are $w_{1},w_{2}\in U$
with $[A(w_{1})]\neq [A(w_{2})]$ such that
$m([A(w_{1})])=m([A(w_{2})])$ for all moments $m\in {\cal M}$.
We set $A_{n}(w):=(a_{1}(w,\cdot),\dots, a_{n}(w,\cdot), 0,\dots)\in X$. 
Then it is easy
to see that the previous condition is equivalent to 
\begin{equation}\label{moment2}
m([A_{n}(w_{1})])=m([A_{n}(w_{2})])\ \ \ {\rm for\ all}\ \ \ n\in\N\ ,\ \ 
m\in {\cal M}\ .
\end{equation}
We will show that $A_{n}(w)\in X(D_{n})$ for all $w\in U$. Then
from (\ref{moment2}) and Theorem \ref{te2} we will get 
$[A_{n}(w_{1})]=[A_{n}(w_{2})]$ for all $n$. 
Clearly this is equivalent to the fact that
$[A(w_{1})]=[A(w_{2})]$ in contradiction to our assumption. This will
complete the proof of the proposition.

So, let us show that $A_{n}(w)\in X(D_{n})$. We consider several
cases.
\begin{itemize}
\item[{\bf (A)}]
All coefficients $a_{i}(w,x)$ in the definition of $A_{n}(w)$
are polynomials in $x$.
\end{itemize}

Replacing $x$ by $z\in\Co$ we obtain holomorphic polynomials
$a_{i}(w,z)$ on $\Co$. We can extend 
definitions (\ref{moment1}) and (\ref{e6}) in this case:
$$
\widetilde a_{i}(w,x):=\int_{0}^{x}a_{i}(w,z)\ \!dz\ \ \ \ {\rm and}
\ \ \ \widetilde A_{A_{n}(w)}(z):=(\widetilde a_{i}(w,z),\dots,\widetilde a_{n}(w,z))\ ,
\ \ \ z\in\Co\ .
$$
Since $\widetilde a_{i}(w,\cdot)$ are holomorphic polynomials,
the map $\widetilde A_{A_{n}(w)}$ is holomorphic and polynomial. 
In particular, $C:=\widetilde A_{A_{n}(w)}(\Co)$ is a (possibly singular) 
rational curve in $\Co^{n}$.
It is easy to check that in this case the normalization $C_{n}$ of $C$ is
$\Co$. Thus conditions $(C_{n1})$ and $(C_{n2})$ are valid in this case,
cf. Remark \ref{r1}.
\begin{itemize}
\item[{\bf (B)}]
All coefficients $a_{i}(w,x)$ in the definition of $A_{n}(w)$ 
are polynomials in
$e^{\pm 2\pi ix/T}$ such that
$\widetilde a_{i}(w,\cdot)$ defined by
(\ref{moment1}) are $T$-periodic functions.
\end{itemize}

In this case as above we can extend $\widetilde A_{A_{n}(w)}$ to a 
holomorphic map $\Co\to\Co^{n}$. The coordinates of this map are polynomials 
in $e^{\pm 2\pi iz/T}$. Let $\Co\to\Co^{*}$, $z\mapsto e^{2\pi iz/T}$, be the 
covering map; here $\Co^{*}$ is the quotient of $\Co$ by the action of the
group $T\cdot\Z$ by translations along the $x$-axis. The map
$\widetilde A_{A_{n}(w)}$ is invariant with respect to this action and 
therefore it determines a holomorphic map 
$A_{A_{n}(w)}':\Co^{*}\to\Co^{n}$ whose 
pullback to $\Co$ coincides with $\widetilde A_{A_{n}(w)}$. By the 
definition the
coordinates of $A_{A_{n}(w)}'$ are Laurent polynomials. In particular,
$A_{A_{n}(w)}'$ is algebraic and the Zariski closure of
the image $A_{A_{n}(w)}'(\Co^{*})$ in $\Co^{n}$ is a (possibly singular)
rational curve $C$. It is easy to see that in this case the normalization 
$C_{n}$ of $C$ is either $\Co$ or $\Co^{*}$. Thus according to
Remark \ref{r1}, $A_{n}(w)\in X(D_{n})$ in this case.
\begin{itemize}
\item[{\bf (C)}]
All coefficients $a_{i}(w,x)$ in the definition of $A_{n}(w)$ 
are polynomials in
$e^{\pm 2\pi ix/T}$ and at least one of $\widetilde a_{i}(w,\cdot)$ in the
definition of $\widetilde A_{A_{n}(w)}$ is not $T$-periodic. 
\end{itemize}

In this case we have
$$
\widetilde a_{i}(w,z)=c_{i}(w)z+b_{i}(w,z)
$$
where $b_{i}(w,\cdot)$ are polynomials in $e^{\pm 2\pi i z/T}$ and at least
one of $c_{i}(w)$ is not zero.
\begin{Lm}\label{le5}
$C=\widetilde A_{A_{n}(w)}(\Co)$ belongs to a closed complex curve in 
$\Co^{n}$.
\end{Lm}
{\bf Proof.} Assume without loss of generality that
$c_{1}(w)\neq 0$. Let us define functions $a_{i}'$ by the formulas
$$
a_{1}'(w,z):=\frac{\widetilde a_{1}(w,z)}{c_{1}(w)}\ \ \ {\rm and}\ \ \ 
a_{i}'(w,z):=\widetilde a_{i}(w,z)-c_{i}(w)\ \!a_{1}'(w,z)\ ,
\ \ \ 2\leq i\leq n\ .
$$
Next, define $A':\Co\to\Co^{n}$ by
$$
A'(z):=(a_{1}'(w,z),\dots, a_{n}'(w,z))\ ,\ \ \ z\in\Co\ .
$$
By the definition the (invertible) linear
transformation $L:\Co^{n}\to\Co^{n}$ given by
$$
L(z):=(c_{1}(w)z_{1},z_{2}+c_{2}(w)z_{1},\dots,
z_{n}+c_{n}(w)z_{1})\ ,\ \ \ z=(z_{1},\dots, z_{n})\in\Co^{n}\ ,
$$
maps $C':=A'(\Co)$ onto $C$. Thus in order to prove that $C$ belongs to a 
closed complex curve in $\Co^{n}$ it suffices to check the same for 
$C'$. So without loss of generality we may assume that 
$\widetilde a_{i}=a_{i}'$, that is,
$$
\widetilde a_{1}(w,z)=z+b_{1}(w,z)\ \ \ {\rm and}\ \ \ 
\widetilde a_{i}(w,z)=b_{i}(w,z)\ ,\ \ \
2\leq i\leq n\ .
$$

Representing $\Co^{n}$ as $\Co\times\Co^{n-1}$ 
we write $\widetilde A_{A_{n}(w)}=(A_{1},A_{2}):\Co\to\Co\times\Co^{n-1}$
where $A_{1}:=\widetilde a_{1}(w,\cdot)$ and $A_{2}:\Co\to\Co^{n-1}$ is a
holomorphic $T$-periodic map. Let $p:\Co^{n}\to\Co^{n-1}$ be the natural 
projection. We will consider three cases.
\begin{itemize}
\item[(1)] 
Assume that at least one of $b_{i}(w,z)$ in the definition of $A_{2}$
has nonzero terms involving $e^{2\pi i dz/T}$ and $e^{-2\pi i mz/T}$
for some positive integers $d$ and $m$.
\end{itemize}

Let $r:\Co\to\Co^{*}$, $r(z):=e^{2\pi iz/T}$, be the covering map, and
$A_{2}':\Co^{*}\to\Co^{n-1}$ be the holomorphic map such that
$A_{2}:=A_{2}'\circ r$. The above assumption implies that $A_{2}'$
is proper. We prove that 
$\widetilde A_{A_{n}(w)}:\Co\to\Co^{n}$ is also proper.

For otherwise, there exists a bounded convergent sequence 
$\{z_{l}\}\subset C\subset\Co^{n}$ such that 
$\widetilde A_{A_{n}(w)}^{-1}(\{z_{l}\})\subset\Co$ is unbounded.
By the definition of $A_{2}$
$$
\widetilde A_{A_{n}(w)}^{-1}(\{z_{l}\})\subset 
A_{2}^{-1}(\{p(z_{l})\}):=\{y_{l}\}+T\cdot\Z\subset\Co
$$
where $\{y_{l}\}$ is a sequence of points in the strip 
$S_{T}=\{z\in\Co\ :\ {\rm Re}\ \! z\in [0,T)\}$. Since $A_{2}'$ is
proper, $\{y_{l}\}\subset S_{T}$ is bounded. We write 
$A_{1}(z)=z+h_{1}(z)+h_{2}(z)$ where
$h_{1}$ is a polynomial in $e^{2\pi i z/T}$ and $h_{2}$ is a polynomial
in $e^{-2\pi i z/T}$. Then $h_{1}+h_{2}$ is bounded on $\{y_{l}\}+T\cdot\Z$.
This and the boundedness of
$A_{1}(\widetilde A_{A_{n}(w)}^{-1}(\{z_{l}\}))$ imply that $z$ is bounded 
on $\widetilde A_{A_{n}(w)}^{-1}(\{z_{l}\})$. Therefore
$\widetilde A_{A_{n}(w)}^{-1}(\{z_{l}\})$ is bounded, a contradiction.

Since $\widetilde A_{A_{n}(w)}:\Co\to\Co^{n}$ is proper, by the Remmert 
proper map theorem (see, e.g., [GH]) 
$C=\widetilde A_{A_{n}(w)}(\Co)\subset\Co^{n}$
is a closed complex subvariety. 
This completes the proof of the lemma in case (1).
\begin{itemize}
\item[(2)] Assume that all $b_{i}(w,z)$ in the definition of $A_{2}$ are
either polynomials in $e^{2\pi i z/T}$ or in $e^{-2\pi i z/T}$ and
$A_{2}$ is not constant.
\end{itemize}

Assume, e.g., that all $b_{i}(w,z)$ are polynomials in 
$e^{2\pi i z/T}$ (the proof in the case $b_{i}(w,z)$ are polynomials in 
$e^{-2\pi i z/T}$ is similar). 
Let $A_{2}':\Co^{*}\to\Co^{n-1}$ be such that $A_{2}=A_{2}'\circ r$. 
It is extended to
a polynomial map $\Co\to\Co^{n-1}$ (denoted also by $A_{2}'$). We
set $o:=A_{2}'(0)\in\Co^{n-1}$ and $M:=A_{2}^{-1}(o)\subset\Co$. 
Then $M:=\{s_{l}\}+T\cdot\Z$ where $\{s_{l}\}\subset S_{T}$ is finite.
Let us consider the complex line $V:=p^{-1}(o)\subset\Co^{n}$. 
Since the map  $A_{2}':\Co\setminus (A_{2}')^{-1}(o)\to
\Co^{n-1}\setminus\{o\}$ is proper, the same argument as
in the proof of case (1) implies that
$\widetilde A_{A_{n}(w)}:\Co\setminus M\to\Co^{n}\setminus V$
is proper. Therefore by the proper map theorem 
$\widetilde A_{A_{n}(w)}(\Co\setminus M)$ is a closed complex curve in
$\Co^{n}\setminus V$. 

Further, $\widetilde A_{A_{n}(w)}$ maps $M$ onto $V\cap C$. Note that
the definition of $M$ and $A_{1}$ imply that $A_{1}$ is unbounded on every 
infinite subset of $M$. Thus $V\cap C\subset\Co^{n}$ is a 
discrete subset. In particular, $V\cap C\subset\Co^{n}$ is a complex analytic
subset and $\widetilde A_{A_{n}(w)}(\Co\setminus M)$ is 
a closed complex curve in $\Co^{n}\setminus (V\cap C)$. Therefore by the Levi 
extension theorem (see, e.g., [GH, Chapter 3]) the closure of
$\widetilde A_{A_{n}(w)}(\Co\setminus M)$ is a closed complex curve in 
$\Co^{n}$. To complete the proof of the lemma in this case observe that
this closure coincides with $C$.
\begin{itemize}
\item[(3)] $A_{2}$ is constant.
\end{itemize}

Let $o:=A_{2}(\Co)\subset\Co^{n-1}$. Then $C$ is an open everywhere
dense subset of the complex line $\{(z,o)\ :\ z\in\Co\}\subset\Co^{n}$. 

The proof of the lemma is complete.\ \ \ \ \ $\Box$
\begin{Lm}\label{le6}
Let $\widetilde C\subset\Co^{n}$ be the closed complex curve containing
$C=\widetilde A_{A_{n}(w)}(\Co)$ from Lemma \ref{le5}. Then its normalization 
$\widetilde C_{n}$ is biholomorphic to $\Co$.
\end{Lm}
{\bf Proof.} Clearly it suffices to prove the lemma for $C$ satisfying 
conditions (1) or (2). In this case we have $C=\widetilde C$. Let
$n:C_{n}\to C$ be the normalization of $C$. Then there is a holomorphic map 
$\widetilde A_{n}:\Co\to C_{n}$ such that 
$\widetilde A_{A_{n}(w)}=n\circ\widetilde A_{n}$. Recall that 
in the proof of Lemma \ref{le5} we established that 
$\widetilde A_{A_{n}(w)}:\Co\setminus\widetilde A_{A_{n}(w)}^{-1}(S)\to 
C\setminus S$ is proper for a discrete set $S\subset C$. This implies easily 
that $\widetilde A_{n}:\Co\to C_{n}$ is a proper surjective map. 
In particular, the
induced homomorphism $\widetilde A_{n}^{*}:H^{1}(C_{n},\Co)\to H^{1}(\Co,\Co)$
is injective. Thus $H^{1}(C_{n},\Co)=0$ and, since $C_{n}$ is non-compact,
$C_{n}\cong\Co$. \ \ \ \ \ $\Box$

Lemmae \ref{le5} and \ref{le6} show that $A_{n}(w)\in X(D)$ in case (C).
\begin{itemize}
\item[{\bf (D)}]
At least one of the coefficients $a_{i}(w,x)$ in the definition of
$A_{n}(w)$ is a non-constant polynomial in $x$.
\end{itemize}

It is easy to see that in this case the extended map 
$\widetilde A_{A_{n}(w)}:\Co\to\Co^{n}$ is proper. Then by the Remmert theorem
$\widetilde A_{A_{n}(w)}(\Co)\subset\Co^{n}$ is a closed complex subvariety.
Similarly to the proof of Lemma \ref{le6} we get that the normalization of
$C$ is $\Co$. 

Thus we have proved that for all possible cases $A_{n}(w)\in X(D_{n})$. 
As we explained
above this implies the required statement of the proposition.\ \ \ \ \ $\Box$

Using this proposition let us complete the proof of the theorem. 

By the definition for every 
$m\in {\cal M}$ its pullback $A^{*}m$ is a holomorphic function on $U$. 
Let $K\subset U$ be a compact subset. Since $U$ is Stein, without loss of
generality we will assume that $K$ is holomorphically convex in $U$. Then
we can choose a relatively compact Stein neighbourhood $V\subset U$ of $K$.
Consider a closed complex subvariety $Z$ of $U\times U$ given by 
equations $A^{*}m(z)-A^{*}m(w)=0$, $z\times w\in U\times U$,
$m\in {\cal M}$. By so-called Cartan's Theorem B, and
by our choice of $V$ we obtain that there is a number $l=l(K)\in\N$ 
such that the set $Z\cap (V\times V)$ coincides with 
$\{z\times w\in V\times V\ :\
A^{*}m_{i}(z)-A^{*}m_{i}(w)=0\ ,\ 1\leq i\leq l\}$. 
Further, according to Remark \ref{r2} the set of universal centers 
${\cal U}_{A}$ of equation (\ref{e3}) satisfying the hypotheses of 
Theorem \ref{te1}
is determined as the set of zeros of holomorphic functions 
$A^{*}m$ on $U$ ($m\in {\cal M}$). As above we can find a number
$s=s(K)\in\N$ such that ${\cal U}_{A}\cap V$ coincides with the set of
zeros of $A^{*}m_{i}$, $1\leq i\leq s$, in $V$. Without loss of generality
we will assume that $l=s$.

Let us consider the holomorphic map 
$$
M_{l}=(A^{*}m_{1},\dots,A^{*}m_{l}): U\to\Co^{l}\ .
$$
By $[\ ]:X\to G(X)$, $a\mapsto [a]$, we denote the canonical 
surjection, see section 2.2. Also, by ${\cal A}_{K}\subset {\cal A}$ we 
denote the set $\{[A(w)]\in G(X)\ :\ w\in K\}$.
\begin{Lm}\label{le7}
There is a continuous map $\widetilde M_{l}:{\cal A}\to\Co^{l}$ which
embeds ${\cal A}_{K}$ into $\Co^{l}$ and satisfies 
$M_{l}=\widetilde M_{l}\circ [\ ]\circ A$. 
\end{Lm}
{\bf Proof.} Since every moment is the pullback with respect to $[\ ]$ of a 
continuous function on $G(X)$, the map $\widetilde M_{l}:{\cal A}\to\Co^{l}$,
$\widetilde M_{l}([a]):=M_{l}(a)$, $a\in A(U)$, is well defined and 
continuous. Let us show that $\widetilde M_{l}$ separates points on 
${\cal A}_{K}$. Let $g_{1},g_{2}\in {\cal A}_{K}$ satisfy
$\widetilde M_{l}(g_{1})=\widetilde M_{l}(g_{2})$. Assume that
$w_{i}\in K$ is such that $[A(w_{i})]=g_{i}$, $i=1,2$. Then by our
assumption $m_{i}(A(w_{1}))=m_{i}(A(w_{2}))$ for all $1\leq i\leq l$.
From here by the choice of $l$ we get $m(A(w_{1}))=m(A(w_{2}))$ for all
$m\in {\cal M}$. Then Proposition \ref{pr1} implies that $g_{1}=g_{2}$.
\ \ \ \ \ $\Box$
\begin{R}\label{comp}
{\rm
Since $K\subset U$ is compact, ${\cal A}_{K}\subset G(X)$ is a 
compact subset by the definition of the topology on $G(X)$, see 
[Br3, Theorem 2.4]. Thus from the lemma we get
$\widetilde M_{l}:{\cal A}_{K}\to M_{l}(K)$ is a homeomorphism.}
\end{R}

Let $c_{n}$, $n\in\N$, be iterated polynomials on $X$ whose set of zeros is
the center set ${\cal C}$, see section 2.1. By the definition, 
the set of zeros
of holomorphic functions $A^{*}c_{n}$ on $U$ is ${\cal C}_{A}$. Next, by
our choice of $V\subset U$, there is a number $N=N(K)\in\N$ such that
${\cal C}_{A}\cap V$ coincides with the set of zeros of the family
$\{A^{*}c_{i}\}_{1\leq i \leq N}$, in $U$. 
Let $\widehat c_{i}$ be a continuous function on
$G(X)$ whose pullback to $X$ coincides with $c_{i}$. According to
Lemma \ref{le7} and the Tietze-Urysohn extension theorem 
applied to the continuous functions 
$(\widetilde M_{l}^{-1})^{*}(\widehat c_{i}|_{{\cal A}_{K}})$ 
on the compact set
$A_{l}(K)\subset\Co^{l}$, there are continuous
functions $f_{i}$ on $\Co^{l}$ such that 
$\widetilde M_{l}^{*}f_{i}|_{{\cal A}_{K}}=\widehat c_{i}|_{{\cal A}_{K}}$,
$1\leq i\leq N$. Finally, the application of Lemma \ref{le7} leads to the
required statement. Namely, under the hypotheses of the theorem the center 
set ${\cal C}_{A}\cap K$ of (\ref{e3}) in
$K$ coincides with the set of zeros of continuous functions 
$f_{n}(A^{*}m_{1},\dots, A^{*}m_{l})$ on $U$ ($1\leq n\leq N$). Also, by
our choice of $l$ the set of universal centers 
${\cal U}_{A}\cap K$ of (\ref{e3}) in $K$
is determined as the set of zeros of the holomorphic functions
$A^{*}m_{n}$, $1\leq n\leq l$, restricted to $K$. 
\ \ \ \ \ $\Box$
\sect{\hspace*{-1em}. Proofs of Theorems \ref{cor1} and \ref{cor2}.}
{\bf 5.1. Proof of Theorem \ref{cor1}.}
Assume that all coefficients $a_{i}$ in (\ref{e3}) are $\alpha$-homogeneous
polynomials in $w\in\Co^{k}$ for some 
$\alpha=(\alpha_{1},\dots,\alpha_{k})\in\N^{k}$ whose coefficients are
either polynomials in $x$ or in $e^{\pm 2\pi ix/T}$. Let us consider
the holomorphic map $P_{\alpha}:\Co^{k}\to\Co^{k}$, 
$$
P_{\alpha}(w_{1},\dots,w_{k}):=(w_{1}^{\alpha_{1}},\dots, 
w_{k}^{\alpha_{k}})\ ,\ \ \ (w_{1},\dots, w_{k})\in\Co^{k}\ .
$$ 
By $\widehat a_{i}(\cdot, x)$ we denote the pullback of $a_{i}(\cdot, x)$
to $\Co^{k}$ with respect to $P_{\alpha}$. Then by the definition all
$\widehat a_{i}$ are homogeneous. Let us consider the equation
(\ref{e3}) whose coefficients are $\widehat a_{i}$. By 
$\widehat A:=(\widehat a_{1}(w,\cdot),\dots)$ we denote the corresponding
map of $\Co^{k}$ into $X$, and by ${\cal U}_{\widehat A}$ and 
${\cal C}_{\widehat A}$ the sets of universal centers and centers for
this equation. Now the hypotheses of the theorem imply that 
${\cal U}_{\widehat A}=\{0\}\subset\Co^{k}$. Also, we clearly have 
\begin{equation}\label{preimage}
{\cal C}_{\widehat A}=P_{\alpha}^{-1}({\cal C}_{A})\ \ \ {\rm and}\ \ \
{\cal U}_{\widehat A}=P_{\alpha}^{-1}({\cal U}_{A})\ .
\end{equation}

Suppose that we will prove that 
${\cal C}_{\widehat A}\subset\Co^{k}$ is determined
as the set of zeros of holomorphic polynomials
$p_{n}(\widehat A^{*}m_{1},\dots, \widehat A^{*}m_{N})$, $1\leq n\leq N$, on
$\Co^{k}$, where $p_{n}$ are some holomorphic polynomials on $\Co^{N}$.
By the definition, see (\ref{moment}), 
$\widehat A^{*}m_{i}$ is the pullback with respect
to $P_{\alpha}$ of the holomorphic polynomial $A^{*}m_{i}$, $1\leq i\leq N$. 
Thus 
$p_{n}(\widehat A^{*}m_{1},\dots, \widehat A^{*}m_{N})$ is the pullback 
with respect to $P_{\alpha}$ of $p_{n}(A^{*}m_{1},\dots, A^{*}m_{N})$, 
$1\leq n\leq N$. From here and (\ref{preimage}) we get that
${\cal C}_{A}$ is determined as the set of zeros of holomorphic 
polynomials $p_{n}(A^{*}m_{1},\dots, A^{*}m_{N})$, $1\leq n\leq N$, on
$\Co^{k}$. This will complete the proof of the theorem. 

Thus without loss of generality we may assume that the coefficients
$a_{i}$ in the original equation (\ref{e3}) are homogeneous polynomials in
$w\in\Co^{k}$.

Further, let us consider a complex algebraic subvariety $Z$ of 
$\Co^{k}\times\Co^{k}$ given by polynomial equations
$A^{*}m(z)-A^{*}m(w)=0$, $z\times w\in\Co^{k}\times\Co^{k}$, $m\in {\cal M}$.
By the Hilbert finiteness theorem there is a number $l\in\N$ such that
$Z=\{z\times v\in\Co^{k}\times\Co^{k}\ :\
A^{*}m_{i}(z)-A^{*}m_{i}(w)=0\ ,\ 1\leq i\leq l\}$. 
Also, according to Remark \ref{r2} the set of universal centers 
${\cal U}_{A}$ of equation (\ref{e3}) satisfying the hypotheses of 
Theorem \ref{te1} is determined as the set of zeros of homogeneous 
holomorphic polynomials
$A^{*}m$ on $\Co^{k}$ ($m\in {\cal M}$). As above we can find a number
$s\in\N$ such that ${\cal U}_{A}$ coincides with the set of
zeros of $A^{*}m_{i}$, $1\leq i\leq s$. Without loss of generality
we will assume that $l=s$.

Let us consider the polynomial map
\begin{equation}\label{poly}
M_{l}=(A^{*}m_{1},\dots,A^{*}m_{l}):\Co^{k}\to\Co^{l}\ .
\end{equation}

We prove
\begin{Lm}\label{proper1}
Under the hypotheses of the theorem $M_{l}$ is a proper map.
\end{Lm}
{\bf Proof.} Assume, to the contrary, that $M_{l}$ is not proper. 
Then there exists an 
unbounded sequence $\{w_{j}\}\subset\Co^{k}$ such that
$\{M_{l}(w_{j})\}$ converges to some $w\in\Co^{l}$. We write
$w_{j}=|w_{j}|\cdot v_{j}$ where $|\cdot|$ is the Euclidean norm on
$\Co^{k}$ and $v_{j}$ belongs to the unit sphere $S^{2k-1}\subset\Co^{k}$.
Without loss of generality we will assume that $\{v_{j}\}$ converges to a
point $v\in S^{2k-1}$. Now,
since every $A^{*}m_{i}$ is a homogeneous polynomial, we have
$A^{*}m_{i}(w_{j})=|w_{j}|^{d_{i}}A^{*}m_{i}(v_{j})$, where
$d_{i}$ is the degree of $A^{*}m_{i}$. According to our assumptions
$\{A^{*}m_{i}(w_{j})\}\subset\Co$ converges and $|w_{j}|\to\infty$ as
$j\to\infty$. This implies easily that $A^{*}m_{i}(v)=0$ for all 
$1\leq i\leq l$. But then by the choice of $l$ we have $A^{*}m(v)=0$ for
all $m\in {\cal M}$. Thus by Remark \ref{r2} we get from here that
$v\in {\cal U}_{A}$ which
contradicts to ${\cal U}_{A}=\{0\}$. This contradiction proves the 
lemma.\ \ \ \ \ $\Box$

Since ${\cal C}_{A}\subset\Co^{k}$ is a complex algebraic subvariety,
from this lemma and the fact that $M_{l}$ is polynomial it follows that
$M_{l}({\cal C}_{A})\subset\Co^{l}$ is a complex algebraic subvariety, as 
well, see, e.g., [GH]. Hence, by the Chow theorem there are holomorphic
polynomials $p_{1},\dots, p_{N}$ on $\Co^{l}$ whose set of zeros is
$M_{l}({\cal C}_{A})$. Then the set of zeros of their pullbacks 
$p_{n}(A^{*}m_{1},\dots, A^{*}m_{N})$ with respect to
$M_{l}$ to $\Co^{k}$ determines the set 
$M_{l}^{-1}(M_{l}({\cal C}_{A}))$. Let us show now that 
\begin{equation}\label{ident}
M_{l}^{-1}(M_{l}({\cal C}_{A}))={\cal C}_{A}\ .
\end{equation}

Assume, to the contrary, that there is 
$w\in M_{l}^{-1}(M_{l}({\cal C}_{A}))\setminus {\cal C}_{A}$. By the
definition there exists some $v\in {\cal C}_{A}$ such that 
$M_{l}(w)=M_{l}(v)$. This implies that $A^{*}m_{i}(w)=A^{*}m_{i}(v)$
for all $1\leq i\leq l$. But then by our choice of $l$ the similar identity
holds for every $m\in {\cal M}$. Now according to Proposition \ref{pr1}
we get from here $[A(w)]=[A(v)]$ in $G(X)$. Then from Corollary \ref{c2}
we obtain that $w\in {\cal C}_{A}$. This contradiction proves (\ref{ident})
and completes the proof of the theorem.\ \ \ \ \ $\Box$\\
{\bf 5.2. Proof of Theorem \ref{cor2}.}
As in the proof of Theorem \ref{cor1} we may assume without loss of
generality that all coefficients $a_{i}$ in (\ref{e3}) are homogeneous
polynomials in $w\in\Co^{k}$ whose coefficients are either polynomials in $x$
or in $e^{\pm 2\pi ix/T}$. Let us consider the polynomial map
$M_{l}:\Co^{k}\to\Co^{l}$ constructed in the same way
as in the proof of Theorem \ref{cor1}, see (\ref{poly}). 
Let $V\subset\Co^{k}$ be
a complex algebraic subvariety whose closure in $\Co\P^{k}$ does not meet
$\overline{\cal U}_{A}$. Then we prove
\begin{Lm}\label{proper2}
Under the hypotheses of the theorem $M_{l}|_{V}$ is a proper map.
\end{Lm}
{\bf Proof.} 
Let $p:\Co^{k}\setminus\{0\}\to\Co\P^{k-1}$, $(w_{1},\dots, w_{k})\mapsto
(w_{1}:\dots : w_{k})$, $(w_{1},\dots w_{k})\in\Co^{k}$, be the canonical
projection determining $\Co\P^{k-1}$. By the definition $p$ can be
extended to a holomorphic map 
$\widehat p:\Co\P^{k}\setminus\{0\}\to\Co\P^{k-1}$ such that
$\widehat p|_{H_{k}}:H_{k}\to\Co\P^{k-1}$ is biholomorphic. (Here
$H_{k}\subset\Co\P^{k}$ is the hyperplane at infinity.) Further, since
${\cal U}_{A}$ is determined as the set of 
zeros of homogeneous polynomials $A^{*}m_{i}$, $1\leq i\leq l$, see 
Remark \ref{r2}, $\widehat p\ \!(\overline{\cal U}_{A}\setminus\{0\})$ is a 
projective variety isomorphic to $\overline{\cal U}_{A}\cap H_{k}$.

Assume, to the contrary, that $M_{l}|_{V}$ is not proper. 
Then there is an 
unbounded sequence $\{w_{j}\}\subset V$ such that
$\{M_{l}(w_{j})\}$ converges in $\Co^{l}$. Without loss of
generality we will assume that $\{w_{j}\}$ converges to a point
$\overline{w}\in H_{k}\cap\overline{V}$. By our assumption
$\overline{w}\not\in \overline{\cal U}_{A}$. We write
$w_{j}=|w_{j}|\cdot v_{j}$ where $|\cdot|$ is the Euclidean norm on
$\Co^{k}$ and $v_{j}$ belongs to the unit sphere $S^{2k-1}\subset\Co^{k}$.
Without loss of generality we will assume that $\{v_{j}\}$ converges to a
point $v\in S^{2k-1}$. Then as in the proof of Lemma \ref{proper1} we
get $v\in {\cal U}_{A}$. Further, by our construction, 
$\widehat p(v)=\widehat p(\overline{w})$. Hence, we have
$\widehat p(\overline{w})\in
\widehat p\ \!(\overline{\cal U}_{A}\setminus\{0\})$.
From here and the facts that  $\widehat p: H_{k}\to\Co\P^{k-1}$ and
$\widehat p: \overline{\cal U}_{A}\cap H_{k}\to  
\widehat p\ \!(\overline{\cal U}_{A}\setminus\{0\})$
are isomorphisms we derive $\overline{w}\in\overline{\cal U}_{A}$.
This contradiction proves the lemma.\ \ \ \ \ $\Box$

Since ${\cal C}_{A}\cap V\subset\Co^{k}$ is a complex algebraic subvariety,
from this Lemma we get that $M_{l}({\cal C}_{A}\cap V)\subset\Co^{l}$
is a complex algebraic subvariety, as well. Thus by the Chow theorem
$M_{l}({\cal C}_{A}\cap V)$ is determined as the set of zeros of
holomorphic polynomials on $\Co^{l}$. Finally, an argument similar to that 
used in the proof of Theorem \ref{cor1} shows that the zero locus of the 
pullbacks of these polynomials to $V$ with respect to $M_{l}$ coincides with
${\cal C}_{A}\cap V$.\ \ \ \ \ $\Box$\\
\\
{\em Acknowledgement.}  The author thanks Y. Yomdin for inspiring
discussions.

\end{document}